%% file: pcompleteDP.tex
\documentclass{article}

\usepackage{proof}
\usepackage{latexsym}
\usepackage{amssymb}

\input{header}

\title{A polynomial time complete disjunction property in intuitionistic propositional logic}
\author{Toshiyasu Arai
\\
Graduate School of Science,
Chiba University
\\
1-33, Yayoi-cho, Inage-ku,
Chiba, 263-8522, JAPAN
\\
tosarai@faculty.chiba-u.jp
}
\date{}
\begin{document}
\maketitle

\begin{abstract}
We extend the polynomial time algorithms due to Buss and Mints\cite{BussMints} and 
 Ferrari, Fiorentini and Fiorino\cite{Ferrari} 
to yield a polynomial time complete disjunction property in intuitionistic propositional logic.
\end{abstract}


The {\it disjunction property, DP\/} of the intuitionistic propositional logic {\sf Ip} says that
if a disjunction $\alp_{0}\lor\alp_{1}$ is derivable intuitionistically, then so is $\alp_{i}$ for an $i$.
This property follows from cut-elimination in sequent calculi, normalization theorem in natural deduction,
Kleene's or Aczel's slash $\Gam |C$ or completeness for Kripke models.

Buss and Mints\cite{BussMints} gave a polynomial time algorithm,
which extracts an $i$ from a given derivation of $\alp_{0}\lor\alp_{1}$ in natural deduction 
such that $\alp_{i}$ is intuitionistically valid.
Such a feasible algorithm based on sequent calculi is given in Buss and Pudl\'ak\cite{BussPudlak},
and Ferrari, Fiorentini and Fiorino\cite{Ferrari} provides an algorithm for derivable sequents $\Gam\Rarw\alp_{0}\lor\alp_{1}$
with sets $\Gam$ of Harrop formulas.

The idea in these algorithms, which comes from \cite{BussMints},
is to prove that one of formulas $\alp_{0}$ and $\alp_{1}$ is in a small set of sequents ({\it immediately derivable sequents\/})
relative to a given intuitionistic derivation of the disjunction $\alp_{0}\lor\alp_{1}$({\bf Boundedness}),
for which there is a polynomial time algorithm testing the membership of sequents in the set,
and any sequent in the set is readily seen to be intuitionistically valid.
In \cite{Ferrari} the authors introduce extraction calculi to generate the set.

In \cite{BussMints} the proof of the {\bf Boundedness} is done
through a partial normalization in natural deduction,
and the proof in \cite{BussPudlak} through cut-elmination.
On the other side, one in \cite{Ferrari} is
 based on an evaluation relation, a variant of Aczel's slash\cite{Aczel}, cf. \cite{MintsKojevnikov}.

In this note we consider the complexity of the DP with Harrop antecedents.
We describe two proofs of {\bf Boundedness}.
One is obtained by a slight modification from \cite{BussMints}, and the other is essentially the same as one in \cite{Ferrari},
but let us stress the fact that the evaluation relation is a feasible restriction of Aczel's slash.
\\

A propositional formula is said to be a {\it Harrop formula\/} if $\lor$ does not occur strictly positive in it.
It is well known that the DP holds with the Harrop antecedents.
Namely for a set of Harrop formulas $\Gam$, if $\Gam\Rarw\alp_{0}\lor\alp_{1}$ is derivable intuitionistically, 
then so is $\Gam\Rarw\alp_{i}$ for an $i$.
Obviously when $\lor$ occurs strictly positive in the antecedent, DP does not hold: $\bet_{0}\lor\bet_{1}\Rarw\alp_{0}\lor\alp_{1}$.
However if we strengthen the antecedent by choosing one $\bet_{i}$ of disjuncts $\bet_{0},\bet_{1}$, 
then one can show $\bet_{i}\Rarw\alp_{n_{i}}$ for some $n_{i}$.
In this way let us generalize the DP:
suppose $\Gam\Rarw\alp_{0}\lor\alp_{1}$ is intuitionistically derivable.
Each {\it strictly positive\/} occurrence of disjunctive subformula $\bet_{0}\lor\bet_{1}$ in the antecedent $\Gam$
 is regarded as an input, and 
we choose $i=0,1$, i.e., one disjunct $\bet_{i}$ freely.
The outputs are disjuncts $\alp_{i}$ from disjunctive formula $\alp_{0}\lor\alp_{1}$ such that
$\Gam^{*}\Rarw\alp_{i}$ is derivable for the strengthening $\Gam^{*}$.
Moreover the problem to find such an $i$ from the given derivation of $\Gam\Rarw\alp_{0}\lor\alp_{1}$
and choices of strictly positive disjuncts in $\Gam$ is solvable in polynomial time, cf. Corollary \ref{th:BM3}.
Suppose that at most one of $\Gam^{*}\Rarw\alp_{i}$ is derivable for any strengthening.
Then the problem has a definite answer.
Indeed, it turn out that the restricted problem is polynomial time complete, cf. Theorem \ref{th:pcomplete}.

\section{Natural deduction {\sf NJp}}

The language of the propositional logic consists of
{\it propositional variables\/} or {\it atoms\/} denoted $p,q,r,\ldots$,
{\it propositional connectives\/} $\bot,\lor,\land,\supset$.
Formulas are denoted by Greek letters $\alp,\bet,\gam,\ldots$
$\lnot\alp:\equiv(\alp\supset\bot)$.
Finite sets of formulas are {\it cedents\/} denoted $\Gam,\Del,\ldots$.
{\it Sequents\/} are ordered pairs of a cedent $\Gam$ and a formula $\alp$ denoted $\Gam\Rarw\alp$,
where $\Gam$ is the {\it antecedent\/} and $\alp$ the {\it succedent formula\/} of the sequent.

\bdf
{\rm Define} strictly positive {\rm occurrence in a (propositional) formula} $\alp$
{\rm in the connectives} $\supset,\lor,\land, \bot$ {\rm as follows.}
{\rm Let} $\bet$ {\rm be an occurrence of a subformula in} $\alp$.
\benu
\item
{\rm If} $\bet\equiv\alp${\rm , then the occurrence is strictly positive in} $\alp$.

\item
{\rm Let} $\alp\equiv(\alp_{0}\lor\alp_{1}), (\alp_{0}\land\alp_{1})$.
 {\rm If the occurrence is strictly positive in} $\alp_{i}${\rm , then the occurrence is strictly positive in} $\alp$.

\item
{\rm Let} $\alp\equiv(\alp_{0}\supset\alp_{1})$.
{\rm If the occurrence is strictly positive in} $\alp_{1}${\rm , then the occurrence is strictly positive in} $\alp$.
\eenu
\edf
The following is the {\bf Axioms} and {\bf Inference rules} in a natural deduction {\sf NJp} for the intuitionistic propositional logic 
{\sf Ip}.
\\

\noindent
{\bf Axioms}.
$\alp,\Gam\Rarw\alp$ for any $\alp$, and $\bot,\Gam\Rarw p$ for atoms $p$.
\\

\noindent
{\bf Inference rules}.
\[
\infer[(\lor E)]{\Gam\Rarw\bet}
{
\Gam\Rarw\alp_{0}\lor\alp_{1}
&
\alp_{0},\Gam\Rarw\bet
&
\alp_{1},\Gam\Rarw\bet
}
\msfiv
\infer[(\lor I)]{\Gam\Rarw\alp_{0}\lor\alp_{1}}
{\Gam\Rarw\alp_{i}}
\]
for $i=0,1$.
\[
\infer[(\land E)]{\Gam\Rarw\alp_{i}}
{\Gam\Rarw\alp_{0}\land\alp_{1}}
\msfiv
\infer[(\land I)]{\Gam\Rarw\alp_{0}\land\alp_{1}}
{
\Gam\Rarw\alp_{0}
&
\Gam\Rarw\alp_{1}
}
\]
\[
\infer[(\supset E)]{\Gam\Rarw\bet}
{
\Gam\Rarw\alp\supset\bet
&
\Gam\Rarw\alp
}
\msfiv
\infer[(\supset I)]{\Gam\Rarw\alp\supset\bet}
{\alp,\Gam\Rarw\bet}
\]
$(\lor I), (\land I), (\supset I)$ are {\it introduction rules\/}, and
$(\lor E), (\land E), (\supset E)$ are {\it elimination rules\/}.

A cedent is {\it Harrop\/} if any formula in it is a Harrop formula.

\section{Immediately derivable sequents}

In this section we introduce the set of {\it immediately derivable sequents from a finite set of sequents\/} $\mathcal{S}$,
and show that there exists a polynomial time algorithm testing whether or not a given sequent is immediately derivable from
$\mathcal{S}$.

\bdf
{\rm Let} $\mathcal{S}$ {\rm be a finite set of sequents. The set of}
immediately derivable {\rm (i.d. $\!$ for short) sequents from} $\mathcal{S}$ {\rm is inductively defined as follows:}
\benu
\item
{\rm Each sequent occurring in} $\mathcal{S}$ {\rm is i.d.$\!$ from} $\mathcal{S}$.
\item
{\rm If both of} $\Gam\Rarw\bet$ {\rm and} $\bet,\Del\Rarw\alp$ {\rm are i.d.$\!$ from} $\mathcal{S}$, 
{\rm then so is}
$\Gam,\Del\Rarw\alp$.
\eenu
\edf

A {\it $(cut)$-deduction\/} is a deduction which may starts with arbitrary sequents and all of whose inference rules are $(cut)$:
\[
\infer[(cut)]{\Gam,\Del\Rarw\alp}
{
\Gam\Rarw\bet
&
\bet,\Del\Rarw\alp
}
\]
Even if we have in hand derivations of $\Gam\Rarw\bet$ and of $\bet,\Del\Rarw\alp$ in {\sf NJp},
the $(cut)$ does not denote the derivation of $\Gam,\Del\Rarw\alp$ obtained by substitution.

Thus a sequent is i.d. $\!$ from $\mathcal{S}$ iff there exists a $(cut)$-deduction stating from sequents in the set $\mathcal{S}$.
\\

A {\it literal\/} is either an atom ({\it positive literal\/}) or its negation({\it negative literal\/}).
A {\it clause\/} is a finite set of literals denoting their disjunction, 
and it is a {\it Horn clause\/} if it contains at most one positive literal.

There is a polynomial time algorithm `HORN SATISFIABILITY' such that given a set $\calh$ of Horn clauses, if it is unsatisfiable, then it returns
a positive unit
 resolution refutation (unit propagation) of $\calh$, and otherwise it returns `SATISFIABLE':
 For a Horn clause $C=\{\ell_{0},\ldots,\ell_{m}\}$ and a positive literal $p$, let
 $C_{p}:=\{\ell_{i}: \ell_{i}\neq\bar{p}, i\leq m\}$.
Define recursively $\calh_{n}$ as follows.
Let $\calh_{0}=\calh$. 
Having defined $\calh_{n}$, if $\calh_{n}$ contains the empty clause $\Box$, then `UNSATISFIABLE'.
Suppose $\Box\not\in\calh_{n}$.
Pick a positive literal $p$ from $\calh_{n}$ if such a $p$ exists,
and let $\calh_{n+1}=\{C_{p}:C\in\calh_{n}\}$.
Otherwise `SATISFIABLE'.
The process $\calh_{n}\leadsto\calh_{n+1}$ is performed at most $N$-times for the number $N$ of atoms
occurring in $\calh$.
Furthermore the cardinality of the sets $\calh_{n}$ of clauses is at most one of $\calh$.
Hence the running time of the algorithm is bounded by square of the size of $\calh$.

\bprp\label{prp:BM1}
There exists a polynomial time algorithm running as follows.
Suppose a finite set $\mathcal{S}$ of sequents and a sequent $\Gam\Rarw\alp$ are given.
If there exists a subsequent $\Gam^{\prime}\Rarw\alp\,(\Gam^{\prime}\subset\Gam)$ of $\Gam\Rarw\alp$ 
which is i.d.$\!$ from $\mathcal{S}$, 
then the algorithm returns a $(cut)$-deduction of a subsequent $\Gam^{\prime}\Rarw\alp$
from $\mathcal{S}$.
Otherwise it returns `NO'.
\eprp
\bprf
Given a finite set $\mathcal{S}$ of sequents and a sequent $\Gam\Rarw\alp$,
consider the set of Horn clauses $\mathcal{S}\cup\{\Rarw \bet: \bet\in\Gam\}\cup\{\alp\Rarw\}$, 
where each formula is regarded as an atom (positive literal).
Then run the algorithm `HORN SATISFIABILITY'.
If the answer is `SATISFIABLE', then any subsequent $\Gam^{\prime}\Rarw\alp$ is not i.d.$\!$ from $\mathcal{S}$.
Return `NO'.
Otherwise
in the refutation, erase the resolution step for $\bet\in\Gam$ and $\alp$:
\[
\infer{\Del\Rarw\gam}
{
\Rarw\bet
&
\bet,\Del\Rarw\gam
}
\leadsto
\bet,\Del\Rarw\gam
\msfiv
\infer{\Rarw}
{\Rarw\alp
&
\alp\Rarw
}
\leadsto
\Rarw\alp
\]
Then it is a $(cut)$-deduction of a subsequent $\Gam^{\prime}\Rarw\alp$ with $\Gam^{\prime}\subset\Gam$
from $\mathcal{S}$.
\eprf

\section{Polynomial time algorithms}

Given a derivation $d$ of a sequent $\Gam\Rarw\alp_{0}\lor\alp_{1}$ with a Harrop antecedent $\Gam$,
the algorithm returns an $i=0,1$ in polynomial time such that $\Gam\Rarw\alp_{i}$ is intuitionistically valid.

\subsection{Normalization with Harrop antecedents}

We follow \cite{BussMints} in this subsection.

\bdf 
{\rm An occurrence of a formula in a derivation is said to be} Harrop maximal
{\rm if it is a conclusion of an introduction rule, and simultaneously a major premiss of an elimination rule
whose lower sequent has Harrop antecedent.}

{\rm Namely in the left of the following figures} $\alp_{0}\lor\alp_{1}$, $\alp_{0}\land\alp_{1}$ {\rm and} $\alp\supset\bet$
{\rm is Harrop maximal if} $\Gam$ {\rm is a Harrop cedent. The right of the figures is the} contratum {\rm of the left.}
\[
\infer[(\lor E)]{\Gam\Rarw\bet}
{
\infer[(\lor I)]{\Gam\Rarw\alp_{0}\lor\alp_{1}}
{\Gam\Rarw\alp_{i}}
&
\alp_{0},\Gam\Rarw\bet
&
\alp_{1},\Gam\Rarw\bet
}
\leadsto
\infer{\Gam\Rarw\bet}
{\Gam\Rarw\alp_{i}
&
\alp_{i},\Gam\Rarw\bet
}
\]
\[
\infer[(\land E)]{\Gam\Rarw\alp_{i}}
{
\infer[(\land I)]{\Gam\Rarw\alp_{0}\land\alp_{1}}
 {
 \Gam\Rarw\alp_{0}
 &
 \Gam\Rarw\alp_{1}
 }
}
\leadsto
\Gam\Rarw\alp_{i}
\]
\[
\infer[(\supset E)]{\Gam\Rarw\bet}
{
\infer[(\supset I)]{\Gam\Rarw\alp\supset\bet}{\alp,\Gam\Rarw\bet}
&
\Gam\Rarw\alp
}
\leadsto
\infer{\Gam\Rarw\bet}
{\Gam\Rarw\alp
&
\alp,\Gam\Rarw\bet
}
\]
{\rm where}
\[
\infer{\Gam\Rarw\bet}
{
\infer*[d_{0}]{\Gam\Rarw\alp}{}
&
\infer*[d_{1}]{\alp,\Gam\Rarw\bet}{}
}
\]
{\rm denotes a natural deduction derivation of} $\Gam\Rarw\bet$
{\rm which is obtained from} $d_{1}$ {\rm by grafting} $d_{0}$ {\rm on each axiom} $\alp,\Gam,\Del\Rarw\alp$
{\rm and deleting} $\alp$ {\rm from antecedents.} 

{\rm If} $\alp\equiv\bot$, {\rm then first convert} $d_{0}$ {\rm to a derivation of} $\Gam\Rarw\gam$ {\rm for an axiom}
$\bot,\Gam,\Del\Rarw\gam$ {\rm in} $d_{1}$.

{\rm A derivation is} Harrop normal {\rm if it contains no Harrop maximal formula (occurrence).}
\edf
{\bf Remark}.
We are not concerned with {\it permutative conversion\/} in normality of natural deduction derivations.

\bdf
{\rm A sequent is} immediately derivable from a derivation {\rm  if it is i.d$\!$ from the set of sequents occurring in the derivation.}
\edf

\bprp\label{prp:BM5}
If a derivation $d$ is contacted by a Harrop maximal formula, then every sequent in the resulting derivation $d^{\prime}$
is i.d.$\!$ from $d$.
\eprp
\bprf
If both of $d_{0}:\Gam\Rarw\alp$ and $d_{1}:\alp,\Del\Rarw\bet$ are subderivations of $d$ (with a Harrop cedent $\Gam$),
then each sequent in a derivation
\[
\infer{\Gam,\Del\Rarw\bet}
{\infer*[d_{0}]{\Gam\Rarw\alp}{}
&
\infer*[d_{1}]{\alp,\Del\Rarw\bet}{}
}
\]
is i.d.$\!$ from $d$.
\eprf

\bprp\label{prp:normalization}
Any {\sf NJp}-derivation $d$ is Harrop normalizable.
Any sequent occurring in the resulting Harrop normal derivation is i.d.$\!$ from $d$.
\eprp
\bprf
This follows from Proposition \ref{prp:BM5}.
\eprf



\bprp\label{prp:BM4}
Let $\Gam$ be a Harrop cedent, $\bot\not\in\Gam$
and $\alp$ a non-Harrop formula (containing a strictly positive occurrence of $\lor$).
Then any Harrop normal derivation of $\Gam\Rarw\alp$ 
 ends with an introduction rule.
\eprp
\bprf
$\Gam\Rarw\alp$ is not an axiom $\alp,\Gam_{0}\Rarw\alp$ nor $\bot,\Gam_{0}\Rarw p$.
Suppose that the derivation ends with an elimination rule.
Consider the left most branch of the derivation tree up to an introduction rule.
Any antecedent on the branch is the Harrop cedent $\Gam$, and any succedent formula $\bet$
contains a strictly positive occurrence of $\lor$.
Therefore there is no introduction rule on the branch. 
 Otherwise the conclusion of the lowest such rule 
would be Harrop maximal.
However the topmost sequent $\Gam\Rarw\bet$ is not an axiom.
This is a contradiction.
\eprf
\\

Suppose $d_{0}$ is an {\sf NJp}-derivation of $\Gam_{0}\Rarw\alp_{0}\lor\alp_{1}$ with a Harrop antecedent $\Gam_{0}$.
If $\bot\in\Gam_{0}$, then pick any $i=I(d)$ and $\Gam_{0}\Rarw\alp_{i}$ is derivable.
Assume $\bot\not\in\Gam_{0}$.

By Proposition \ref{prp:normalization} Harrop normalize the derivation $d_{0}$ to get
 a Harrop normal derivation $d_{1}$ of $\Gam_{0}\Rarw\alp_{0}\lor\alp_{1}$.
Any sequent occurring in $d_{1}$ is i.d.$\!$ from $d_{0}$.
On the other side by Proposition \ref{prp:BM4} $d_{1}$ ends with an introduction rule, i.e., a $(\lor I)$:
\[
\infer[(\lor I)]{\Gam_{0}\Rarw\alp_{0}\lor\alp_{1}}
{\Gam_{0}\Rarw\alp_{i}}
\]
Therefore one of the sequents $\Gam_{0}\Rarw\alp_{i}$ is i.d.$\!$ from $d_{0}$.

Now first check whether or not a subsequent $\Gam_{0}^{\prime}\Rarw\alp_{0}$ is i.d.$\!$ from $d_{0}$ 
by the polynomial time algorithm
in Proposition \ref{prp:BM1}.
If it is the case, then $I(d)=0$.
Otherwise $I(d)=1$.
Moreover in each case the algorithm yields a $(cut)$-deduction of $\Gam^{\prime}\Rarw\alp_{I(d)}$
from the set of sequents occurring in $d_{0}$.

Thus we have shown the
\bth\label{th:BM3H}
There exists a polynomial time algorithm running as follows.
Given an {\sf NJp}-derivation $d$ of a sequent $\Gam\Rarw\alp_{0}\lor\alp_{1}$ with a Harrop antecedent $\Gam$,
the algorithm yields an $i=0,1$ and a $(cut)$-deduction of a subsequent $\Gam^{\prime}\Rarw\alp_{i}$ of 
$\Gam\Rarw\alp_{i}$ 
from the set of sequents occurring in $d$.
\end{theorem}

\subsection{A feasible Aczel's slash}

We follow \cite{Ferrari} in this subsection.

\bdf
{\rm Let} $\mathcal{S}$ {\rm be a set of sequents,} $\Gam$ {\rm a cedent and} $\alp$ {\rm a formula. Then}
$\mathcal{S}:\Gam | \alp$ {\rm holds iff} $\Gam\Rarw\alp$ {\rm is i.d.$\!$ from} $\mathcal{S}${\rm , and one of the following
conditions holds inductively:}
\benu
\item
$\alp\in Var\cup\{\bot\}$.

\item
$\alp\equiv(\bet\supset\gam)$, {\rm and if} 
$\mathcal{S}:\Gam |\bet$
{\rm , then} $\mathcal{S}:\Gam |\gam$.

\item
$\alp\equiv(\bet_{0}\land\bet_{1})$ {\rm and} $\mathcal{S}:\Gam |\bet_{i}$ {\rm for any} $i$.

\item
$\alp\equiv(\bet_{0}\lor\bet_{1})$ {\rm and} $\mathcal{S}:\Gam |\bet_{i}$ {\rm for some} $i$.
\eenu
{\rm For a cedent} $\Del$,
$\mathcal{S}:\Gam |\Del$ {\rm iff} $\mathcal{S}:\Gam |\alp$ {\rm for any} $\alp\in\Del$.
\edf

\bprp\label{prp:soundH1}
{\rm (Soundness)}
Suppose a set $\mathcal{S}$ of sequents contains any sequents occurring in a derivation $d$.
For any sequent $\Gam\Rarw\alp$ occurring in $d$,
if $\mathcal{S}:\Gam_{0} |\Gam$, then $\mathcal{S}:\Gam_{0} |\alp$.
\eprp
\bprf
By induction on the length of subproof of the sequent in $d$.

If $\Gam\Rarw\alp$ occurs in $d$ and $\mathcal{S}:\Gam_{0} |\Gam$, then all of $\Gam\Rarw\alp$ and 
$\Gam_{0}\Rarw\gam$ for $\gam\in\Gam$ are i.d. $\!\!$from $\mathcal{S}$.
Hence so is $\Gam_{0}\Rarw\alp$.
This shows the case when $\Gam\Rarw\alp$ is an axiom $\bot,\Del\Rarw p$.

First consider
\[
\infer[(\lor E)]{\Gam\Rarw\bet}
{
\Gam\Rarw\alp_{0}\lor\alp_{1}
&
\alp_{0},\Gam\Rarw\bet
&
\alp_{1},\Gam\Rarw\bet
}
\]
By IH $\mathcal{S}:\Gam_{0}|(\alp_{0}\lor\alp_{1})$.
Let $i=0,1$ be such that $\mathcal{S}:\Gam_{0}|\alp_{i}$.
IH yields $\mathcal{S}:\Gam_{0}|\bet$.
Next consider
\[
\infer[(\supset I)]{\Gam\Rarw\alp\supset\bet}
{\alp,\Gam\Rarw\bet}
\]
Suppose $\mathcal{S}:\Gam_{0}|\Gam$.
It suffices to show that $\mathcal{S}:\Gam_{0}|\bet$ assuming $\mathcal{S}:\Gam_{0}|\alp$, which follows from IH.

Third consider
\[
\infer[(\supset E)]{\Gam\Rarw\bet}
{
\Gam\Rarw\alp\supset\bet
&
\Gam\Rarw\alp
}
\]
$\Gam\Rarw\bet$ is i.d. $\!\!$from $\mathcal{S}$. 
IH yields $\mathcal{S}:\Gam_{0}|(\alp\supset\bet)$ and $\mathcal{S}:\Gam_{0}|\alp$, and hence
$\mathcal{S}:\Gam_{0}|\bet$.

Other cases are seen easily.
\eprf

\bdf
{\rm For formulas} $\gam$, {\rm a class of sequents} $\mathcal{C}(\gam)$ {\rm is defined recursively.}

\benu
\item
$\mathcal{C}(\alp)=\mathcal{C}(\alp_{0}\lor\alp_{1})=\emptyset$
{\rm for any atomic formula} $\alp\in Var\cup\{\bot\}$, {\rm and any disjunctive formula} $\alp_{0}\lor\alp_{1}$.

\item
$\mathcal{C}(\alp\supset\bet)=\{\alp,\alp\supset\bet\Rarw\bet\}\cup\mathcal{C}(\bet)$.

\item
$\mathcal{C}(\alp_{0}\land\alp_{1})=\{\alp_{0}\land\alp_{1}\Rarw\alp_{i}:i=0,1\}\cup\bigcup_{i=0,1}\mathcal{C}(\alp_{i})$.

\eenu

\edf
It is easy to see that the size of the set $\mathcal{C}(\gam)$ is bounded by a polynomial of the size of $\gam$,
and $\mathcal{C}(\gam)$ is polynomial time recognizable.

\bprp\label{prp:soundH2}
For any Harrop formula $\alp_{0}$ with $\mathcal{C}(\alp_{0})\subset\mathcal{S}$,  
if $\Gam_{0}\Rarw\alp_{0}$ is i.d.$\!$ from $\mathcal{S}$,
then $\mathcal{S}:\Gam_{0} |\alp_{0}$.
\eprp
\bprf
Let $\alp_{0}$ be a Harrop formula such that $\mathcal{C}(\alp_{0})\subset\mathcal{S}$.
We claim for any strictly positive subformula $\alp$ of $\alp_{0}$,
if $\Gam_{0}\Rarw\alp$ is i.d. $\!\!$from $\mathcal{S}$, then $\mathcal{S}:\Gam_{0}|\alp$.


The claim is shown by induction on $\alp$.
There is nothing to show when $\alp$ is an atomic formula in $Var\cup\{\bot\}$.

Let $\alp\equiv(\bet\supset\gam)$ be a formula not of the form (\ref{eq:itesup}).
Since then, $\alp,\bet\Rarw\gam$ is i.d. $\!\!$from $\mathcal{S}$, so is $\bet,\Gam_{0}\Rarw\gam$.
To show $\mathcal{S}:\Gam_{0}|\alp$, assume $\mathcal{S}:\Gam_{0}|\bet$.
Then $\Gam_{0}\Rarw\bet$ is i.d. $\!\!$from $\mathcal{S}$, hence so is $\Gam_{0}\Rarw\gam$.
IH yields $\mathcal{S}:\Gam_{0}|\gam$.

Next let $\alp\equiv(\bet_{0}\land\bet_{1})$.
Since $\alp\Rarw\bet_{i}$ is i.d. $\!\!$from $\mathcal{S}$, so is $\Gam_{0}\Rarw\bet_{i}$.
IH yields $\mathcal{S}:\Gam_{0}|\bet_{i}$ for any $i$, and hence $\mathcal{S}:\Gam_{0}|\alp$.

Since $\alp_{0}$ is a Harrop formula, $\alp\not\equiv(\bet\lor\gam)$.
\eprf
\\

Let $d$ be a derivation of a sequent $\Gam_{0}\Rarw\alp$.
A sequent is said to be {\it immediately derivable with analyses\/} (i.d.a. $\!\!$for short) from $d$
iff the sequent is immediately derivable from sequents occurring in $d$, sequents $\Gam_{0}\Rarw\alp_{0}$ 
for each
$\alp_{0}\in\Gam_{0}$ 
and sequents in $\bigcup\{\mathcal{C}(\gam):\gam\in\Gam_{0}\}$.

Let $d$ be a derivation of a sequent $\Gam_{0}\Rarw\alp_{0}\lor\alp_{1}$ with a Harrop cedent $\Gam_{0}$.
Let $\mathcal{S}_{a}$ be the set of sequents occurring in $d$, sequents $\Gam_{0}\Rarw\alp_{0}$ 
for each
$\alp_{0}\in\Gam_{0}$ 
and sequents in $\bigcup\{\mathcal{C}(\gam):\gam\in\Gam_{0}\}$.

By Proposition \ref{prp:soundH2} we have $\mathcal{S}_{a}:\Gam_{0}|\Gam_{0}$.
Hence by Proposition \ref{prp:soundH1} we obtain $\mathcal{S}_{a}:\Gam_{0}|(\alp_{0}\lor\alp_{1})$.
Let $i=0,1$ be such that $\mathcal{S}_{a}:\Gam_{0}|\alp_{i}$.
Then $\Gam_{0}\Rarw\alp_{i}$ is i.d.a. $\!\!$from $d$.

Check whether or not a subsequent $\Gam_{0}^{\prime}\Rarw\alp_{0}$ is i.d.a.$\!$ from $d$ 
by the polynomial time algorithm
in Proposition \ref{prp:BM1}.
We have shown the following Theorem \ref{th:BM3A} which is slightly weaker than Theorem \ref{th:BM3H}.

\bth\label{th:BM3A}
There exists a polynomial time algorithm running as follows.
Given an {\sf NJp}-derivation $d$ of a sequent $\Gam\Rarw\alp_{0}\lor\alp_{1}$ with a Harrop antecedent $\Gam$,
the algorithm yields an $i=0,1$ and a $(cut)$-deduction of a subsequent $\Gam^{\prime}\Rarw\alp_{i}$ of 
$\Gam\Rarw\alp_{i}$ 
from the set of sequents occurring in $d$, sequents $\Gam\Rarw\gam$ for each
$\gam\in\Gam$ 
and sequents in $\bigcup\{\mathcal{C}(\gam):\gam\in\Gam\}$.
\end{theorem}

\section{A generalized DP}

For a formula $\alp$, $\alp_{\lor}^{+}$ denotes the set of strictly positive 
occurrences of disjunctive formulas in $\alp$.
For a cedent $\Gam$, $\Gam_{\lor}^{+}$ denotes the set of strictly positive 
occurrences of disjunctive formulas in one of formulas in $\Gam$.

$n_{\Gam}$ denote the cardinality of the set $\Gam_{\lor}^{+}$.

Given a cedent $\Gam$,
enumerate the elements in the set $\Gam_{\lor}^{+}$,
$\{\bet^{j}_{0}\lor\bet^{j}_{1}: j<n_{\Gam}\}$ such that if $j_{0}<j_{1}$, then 
$\bet^{j_{0}}_{0}\lor\bet^{j_{0}}_{1}$ is not a subformula of $\bet^{j_{1}}_{0}\lor\bet^{j_{1}}_{1}$.
Each binary number $k=\sum_{j<n_{\Gam}} k_{j}2^{j}<2^{n_{\Gam}}$ is identified with the choice such that
the disjunct $\bet^{j}_{k_{j}}$ is chosen from the disjunction $\bet^{j}_{0}\lor\bet^{j}_{1}$.

Let $\alp(k)$ denote the formula obtained from $\alp\in\Gam$
by replacing the disjunction $\bet^{j}_{0}\lor\bet^{j}_{1}$ by the disjunct $\bet^{j}_{k_{j}}$,
where replacements are done longer subformula occurrences first.
Namely first replace $\bet^{0}_{0}\lor\bet^{0}_{1}$ by $\bet^{0}_{k_{0}}$,
and then $\bet^{1}_{0}\lor\bet^{1}_{1}$ by $\bet^{1}_{k_{1}}$, and so forth.

\bprp\label{prp:strengthen}
$\alp(k)\Rarw\alp$ 
is intuitionistically valid.
\eprp
Suppose $d$ is an {\sf NJp}-derivation of $\Gam\Rarw\alp_{0}\lor\alp_{1}$ and let $k<2^{n_{\Gam}}$.
Then ${\sf NJp}\vdash\Gam(k)\Rarw\alp_{0}\lor\alp_{1}$ by an 
{\sf NJp}-derivation $d(k)$ of $\Gam(k)\Rarw\alp_{0}\lor\alp_{1}$, where $d(k)$ is polynomial time 
computable from $d$ and $k$.
Theorem \ref{th:BM3H} yields the
\bcor\label{th:BM3}
There exists a polynomial time algorithm running as follows.
Given an {\sf NJp}-derivation $d$ of a sequent $\Gam\Rarw\alp_{0}\lor\alp_{1}$ and a number $k<2^{n_{\Gam}}$,
the algorithm yields an $i=0,1$ and a $(cut)$-deduction of a subsequent $\Gam(k)^{\prime}\Rarw\alp_{i}$ of 
$\Gam(k)\Rarw\alp_{i}$ 
from the set of sequents occurring in $d(k)$.
\ecor

\section{Polynomial time completeness}

It is well known that `UNIT' is polynomial time complete 
where `UNIT' is a problem to determine whether or not there is a unit resolution refutation of a given 
set of clauses.
Let us modify the proof of the completeness in \cite{JonesLaaser} to show 
the polynomial time completeness of a generalized DP.

Let $M$ be a deterministic one-tape Turing machine which operates in at most a polynomial $\ell=p(n)$
for inputs of length $n$.
Suppose that $M$ has initial state $s_{0}$, accepting state $s_{a}$ and rejecting state $s_{r}$ such that $M$ eventually reaches one of states $s_{a},s_{r}$, and remains in that state without terminating, scanning a blank $B$ at its starting position.
$M$ never moves to the left of its starting position.
Let $\Sig$, $\Gam$ and $Q$ be a set of input symbols, tape symbols and states, resp.
An {\it instantaneous description, ID\/} is a string $\sig$ on $\Gam\cup(Q\times\Gam)$ in which
symbols in $Q\times\Gam$ occurs exactly once.

Let $P^{a}_{i,t}$ be atoms for $a\in\Gam\cup(Q\times\Gam)$, $i\leq\ell+1$ and $t\leq\ell$.
We write $P(a,i,t)$ for $P^{a}_{i,t}$.
$P(a,i,t)$ is intended to express that `$a$ is the $i$-th symbol of a $t$-th $M$-computation $\sig_{t}$',
where the starting position is $1$.
Define formulas $\bet$, $\del_{i}\, (i\leq \ell+1)$ and $\gam$ as follows.
\[
\bet \equiv  \bigwedge_{t}(P(B,0,t)\land P(B,\ell+1,t))
\]
and
\[
\del_{0}  \equiv 
 \oplus_{a_{1}\in\Sig}P((s_{0},a_{1}),1,0)\land \bigwedge_{1<i\leq n}\oplus_{a\in\Sig}P(a,i,0)\land \bigwedge_{n<i\leq\ell}P(B,i,0)
\]
where $\oplus$ denotes the `excluded or',
$\oplus_{i<n}p_{i}:\equiv\bigvee_{i<n}[p_{i}\land\bigwedge_{j\neq i}\lnot p_{j}]$.


$\del_{0}$ states that an initial configuration is given, and $\bet$ says that positions $0$ and $\ell+1$ are 
always blank in computations.

For an input $x=a_{1}\cdots a_{n}$, let $\del_{0}(x)$ be the formula stating the initial configuration on $x$:
\[
\del_{0}(x) \equiv  P((s_{0},a_{1}),1,0)\land \bigwedge_{1<i\leq n}P(a_{i},i,0)\land \bigwedge_{n<i\leq\ell}P(B,i,0)
.\]

Let $f:(\Gam\cup(Q\times\Gam))^{3}\to\Gam\cup(Q\times\Gam)$ be a function describing the transition function of $M$
as follows.
Assume $P(a,i-1,t)\land P(b,i,t)\land P(c,i+1,t)$.
Then $P(f(a,b,c),i,t+1)$ holds.
For $0\leq t<\ell$, let
\[
\del_{t+1}\equiv\bigwedge_{i}\bigwedge_{a,b,c}[P(a,i-1,t)\land P(b,i,t)\land P(c,i+1,t) \supset P(f(a,b,c),i,t+1)]
.\]

$\del_{0}(x)$ as well as $\del_{t}$ for $t>0$ is a conjunction of Horn clauses.
Let
\beqnarrs
\Gam & \equiv & (\bet\land\bigwedge_{0\leq t\leq \ell}\del_{t})
\\
\Gam(x,t) & \equiv & (\bet\land\del_{0}(x)\land\bigwedge_{0< s\leq t}\del_{s})
\eeqnarrs
Each $\Gam(x,t)$ is satisfiable formula for any $x$ and $t\leq\ell$.

\bprp\label{prp:JL7}
Let $\sig_{0}\vdash\sig_{1}\vdash\cdots\vdash\sig_{\ell}$ be the $M$-computation on an input $x=a_{1}\cdots a_{n}$.
For each $a\in\Gam\cup(Q\times\Gam)$, $i\leq\ell+1$ and $t\leq \ell$,
$a$ is the $i$-th symbol of $\sig_{t}$ iff
$\Gam(x,t)\Rarw P(a,i,t)$ is intuitionistically derivable.
\eprp
\bprf
$\vdash\alp$ means the intuitionistic derivability of $\alp$.
By induction on $t$ we show if $a$ is the $i$-th symbol of $\sig_{t}$, then
$\vdash\Gam(x,t)\Rarw P(a,i,t)$.
The converse is seen from the (classical) soundness of the derivability relation $\vdash$.

The case $i=0$ is trivial.
Suppose the proposition holds for $t$, and $a,b,c$ are the symbols at position $i-1,i,i+1$ in $\sig_{t}$.
Then by IH we have $\vdash\Gam(x,t)\Rarw P(a,i-1,t)\land P(b,i,t)\land P(c,i+1,t)$.
By $\del_{t+1}$ we obtain $\vdash\Gam(x,t+1)\Rarw P(f(a,b,c),i,t+1)$.
\eprf

\bcor\label{cor:JL8}
$M$ accepts an input $x=a_{1}\cdots a_{n}$ iff
$\Gam(x,\ell)
\Rarw P((s_{a},B),1,\ell)$ is intuitionistically derivable.
\ecor

We see that $\Gam \Rarw P((s_{a},B),1,\ell)\lor P((s_{r},B),1,\ell)$ is intuitionistically derivable
from Proposition \ref{prp:JL7}.
Though the size of $\Gam$ is polynomial in $n$, the size of the 
above proof of $\Gam \Rarw P((s_{a},B),1,\ell)\lor P((s_{r},B),1,\ell)$
 is exponential since the proof is based on case distinctions and there are exponentially many inputs.

Let
\[
\gam\equiv 
\bigwedge_{a\neq(s_{a},B),(s_{r},B)}\lnot P(a,1,\ell)
.\]

\bprp
There exists an intuitionistic derivation of $\Gam\land\gam \Rarw P((s_{a},B),1,\ell)\lor P((s_{r},B),1,\ell)$
in size polynomial of $n$.
\eprp
\bprf
Let
\[
\alp_{t}\equiv \bigwedge_{1\leq i\leq \ell}\bigvee_{a\in\Gam\cup(Q\times\Gam)}P(a,i,t)
.\]
We show by induction on $t\leq\ell$ that
\[
\vdash\bigwedge_{0\leq s\leq t}\del_{s}\Rarw\alp_{t}
.\]
We have $\vdash\del_{0}\Rarw\alp_{0}$.
Suppose $\vdash\bigwedge_{0\leq s\leq t}\del_{s}\Rarw\alp_{t}$ for $t<\ell$.
Then by $\del_{t+1}$ we have $\vdash\bigwedge_{0\leq s\leq t+1}\del_{s}\Rarw\alp_{t+1}$.

Since $\ell$ is polynomial in $n$, so is each derivation of $\bigwedge_{0\leq s\leq t}\del_{s}\Rarw\alp_{t}$.
Hence a polynomial size derivation of $\Gam\Rarw\alp_{\ell}$ is obtained, and hence one of
$\Gam\Rarw \bigvee_{a\in\Gam\cup(Q\times\Gam)}P(a,1,\ell)$.
Then by $\gam$, $\Gam\land\gam\Rarw P((s_{a},B),1,\ell)\lor P((s_{r},B),1,\ell)$ has a polysize derivation.
\eprf
\\

Let $\alp\equiv \lnot(P((s_{a},B),1,\ell)\land P((s_{r},B),1,\ell))$, and $\Del=\Gam\cup\{\gam,\alp\}$.

In the formulas in $\Del$, $\bet,\gam,\alp$ as well as $\del_{t}$ for $t>0$ are Harrop formulas.
Strictly positive disjunctions occur only in $\del_{0}$, and 
any choice of one disjunct for each strictly positive disjunctions
yields $\Del(x)=\Gam(x,\ell)\land\gam\land\alp$.

There are intuitionistic derivations $d$ of $\Del\Rarw P((s_{a},B),1,\ell)\lor P((s_{r},B),1,\ell)$ and 
of $\Del\Rarw \lnot(P((s_{a},B),1,\ell)\land P((s_{r},B),1,\ell))$, both of which is of polysize in $n$.
Let $I(x,d)=a,r$ be such that $I(x,d)=a$ iff $\vdash\Del\Rarw P((s_{a},B),1,\ell)$,
and $I(x,d)=r$ iff $\vdash\Del\Rarw P((s_{r},B),1,\ell)$.
By Corollary \ref{th:BM3} the predicate $I$ is polynomial time computable.
On the other side by Corollary \ref{cor:JL8},
$M$ accepts an input $x$ iff
$I(x,d)=a$.
Therefore the predicate $I$ is polynomial time complete.
Thus we have shown the

\bth\label{th:pcomplete}
Let $\Gam$ be a cedent, $\alp_{0},\alp_{1}$ formulas such that
$\Gam(x)\cup\{\lnot(\alp_{0}\land\alp_{1})\}$ is (classically) satisfiable for any strengthening $\Gam(x)$ of $\Gam$
by choosing one disjunct from each strictly positive disjunction in $\Gam$.
Then the problem deciding the $i=0,1$ such that
$\Gam(x)\Rarw\alp_{i}$ is intuitionistically derivable 
from given derivation of $\Gam\Rarw\alp_{0}\lor\alp_{1}$ and $x$
is polynomial time complete.
\end{theorem}

\end{document}

%% file: header.tex
\newcommand{\blem}{\begin{lemma}}
\newcommand{\elem}{\end{lemma}}
\newcommand{\bth}{\begin{theorem}}
\newcommand{\benu}{\begin{enumerate}}
\newcommand{\eenu}{\end{enumerate}}
\newcommand{\bdes}{\begin{description}}
\newcommand{\edes}{\end{description}}
\newcommand{\bdf}{\begin{definition}}
\newcommand{\edf}{\end{definition}}
\newcommand{\bcor}{\begin{cor}}
\newcommand{\ecor}{\end{cor}}
\newcommand{\bprp}{\begin{proposition}}
\newcommand{\eprp}{\end{proposition}}
\newcommand{\bmlem}{\begin{mlemma}}
\newcommand{\emlem}{\end{mlemma}}
\newcommand{\bclm}{\begin{claim}}
\newcommand{\eclm}{\end{claim}}
\newcommand{\bprf}{{\bf Proof}.\hspace{2mm}}

\newcommand{\eprf}{\hspace*{\fill} $\Box$}

\newcommand{\beqn}{\begin{equation}}
\newcommand{\eeqn}{\end{equation}}
\newcommand{\beqnarr}{\begin{eqnarray}}
\newcommand{\eeqnarr}{\end{eqnarray}}
\newcommand{\beqnarrs}{\begin{eqnarray*}}
\newcommand{\eeqnarrs}{\end{eqnarray*}}

\newtheorem{theorem}{Theorem}[section]
\newtheorem{definition}[theorem]{Definition}
\newtheorem{proposition}[theorem]{Proposition}
\newtheorem{lemma}[theorem]{Lemma}
\newtheorem{cor}[theorem]{Corollary}
\newtheorem{mlemma}[theorem]{Main Lemma}
\newtheorem{claim}[theorem]{Claim}

\newcommand{\alp}{\alpha}

\newcommand{\del}{\delta}
\newcommand{\Del}{\Delta}

\newcommand{\bet}{\beta}
\newcommand{\gam}{\gamma}
\newcommand{\Gam}{\Gamma}

\newcommand{\sig}{\sigma}
\newcommand{\Sig}{\Sigma}

\newcommand{\Rarw }{\Rightarrow}

\newcommand{\calh}{{\cal H}}

\newcommand{\msfiv}{\mbox{\hspace{5mm}}}